\documentclass[11pt]{amsart}
\usepackage{amsmath,amsthm,amsfonts,amssymb,latexsym,epsfig}

\newtheorem{theorem}{Theorem}[section]

\newtheorem{cor}{Corollary}[section]

\numberwithin{equation}{section}

\theoremstyle{definition}

\theoremstyle{remark}

\begin{document}
\title{A note on $l^p$ norms of weighted mean matrices}
\author{Peng Gao}
\address{Division of Mathematical Sciences, School of Physical and Mathematical Sciences,
Nanyang Technological University, 637371 Singapore}
\email{penggao@ntu.edu.sg}
\date{August 25, 2008.}


\begin{abstract}
  We present some results concerning the $l^p$ norms of weighted mean matrices. These results can
  be regarded as analogues to a result of Bennett concerning
  weighted Carleman's inequalities.
\end{abstract}

\maketitle
\section{Introduction}
\label{sec 1} \setcounter{equation}{0}

  Suppose throughout that $p\neq 0, \frac{1}{p}+\frac{1}{q}=1$.
  For $p \geq 1$, let $l^p$ be the Banach space of all complex sequences ${\bf a}=(a_n)_{n \geq 1}$ with norm
\begin{equation*}
   ||{\bf a}||_p: =(\sum_{n=1}^{\infty}|a_n|^p)^{1/p} < \infty.
\end{equation*}
  The celebrated
   Hardy's inequality (\cite[Theorem 326]{HLP}) asserts that for $p>1$,
\begin{equation}
\label{eq:1} \sum^{\infty}_{n=1}\Big{|}\frac {1}{n}
\sum^n_{k=1}a_k\Big{|}^p \leq \Big (\frac
{p}{p-1} \Big )^p\sum^\infty_{n=1}|a_n|^p.
\end{equation}
   Hardy's inequality can be regarded as a special case of the
   following inequality:
\begin{equation}
\label{01}
   \Big | \Big |C \cdot {\bf a}\Big | \Big |^p_p =\sum^{\infty}_{n=1} \Big{|}\sum^{\infty}_{k=1}c_{n,k}a_k
    \Big{|}^p \leq U \sum^{\infty}_{n=1}|a_n|^p,
\end{equation}
   in which $C=(c_{n,k})$ and the parameter $p>1$ are assumed
   fixed, and the estimate is to hold for all complex
   sequences ${\bf a} \in l^p$. The $l^{p}$ operator norm of $C$ is
   then defined as
\begin{equation*}
\label{02}
    ||C||_{p,p}=\sup_{||{\bf a}||_p = 1}\Big | \Big |C \cdot {\bf a}\Big | \Big |_p.
\end{equation*}
   It follows that inequality \eqref{01} holds for any ${\bf a} \in l^p$ when $U^{1/p} \geq ||C||_{p,p}$ and fails to hold for some ${\bf a} \in l^p$ when $U^{1/p} <||C||_{p,p}$.
    Hardy's inequality thus asserts that the Ces\'aro matrix
    operator $C$, given by $c_{n,k}=1/n , k\leq n$ and $0$
    otherwise, is bounded on {\it $l^p$} and has norm $\leq
    p/(p-1)$. (The norm is in fact $p/(p-1)$.)

    We say a matrix $A=(a_{n,k})$ is a lower triangular matrix if $a_{n,k}=0$ for $n<k$ and a lower triangular matrix $A$ is a summability matrix if
    $a_{n,k} \geq 0$ and
    $\sum^n_{k=1}a_{n,k}=1$. We say a summability matrix $A$ is a weighted
    mean matrix if its entries satisfy:
\begin{equation}
\label{021}
    a_{n,k}=\lambda_k/\Lambda_n,  ~~ 1 \leq k \leq
    n; \hspace{0.1in} \Lambda_n=\sum^n_{i=1}\lambda_i, \lambda_i \geq 0, \lambda_1>0.
\end{equation}

    Hardy's inequality \eqref{eq:1} now motivates one to
    determine the $l^{p}$ operator norm of an arbitrary summability or weighted mean matrix $A$.
    In \cite{G5}, the author proved the following result:
\begin{theorem}
\label{thm03}
    Let $1<p<\infty$ be fixed. Let $A$ be a weighted mean matrix given by
    \eqref{021}. If for any integer $n \geq 1$, there exists a positive constant
    $0<L<p$ such that
\begin{equation}
\label{024}
    \frac {\Lambda_{n+1}}{\lambda_{n+1}} \leq \frac
    {\Lambda_n}{\lambda_n}  \Big (1- \frac
    {L\lambda_n}{p\Lambda_n} \Big )^{1-p}+\frac {L}{p}~~,
\end{equation}
    then
    $||A||_{p,p} \leq p/(p-L)$.
\end{theorem}

   It is easy to see that the above result implies the following well-known result of Cartlidge \cite{Car} (see also \cite[p. 416, Theorem C]{B1}):
\begin{theorem}
\label{thm02}
    Let $1<p<\infty$ be fixed. Let $A$ be a weighted mean matrix given by
    \eqref{021}. If
\begin{equation}
\label{022}
    L=\sup_n\Big(\frac {\Lambda_{n+1}}{\lambda_{n+1}}-\frac
    {\Lambda_n}{\lambda_n}\Big) < p ~~,
\end{equation}
    then
    $||A||_{p,p} \leq p/(p-L)$.
\end{theorem}

   The above result of Cartlidge are very handy to use when determining
   $l^p$ norms of certain weighted mean matrices.
   We refer the readers to the articles \cite{G}, \cite{Be1}, \cite{G6} and \cite{G5} for more recent
   developments in this area.

We note here that by a change of variables $a_k \rightarrow
a^{1/p}_k$ in \eqref{eq:1} and on letting $p \rightarrow +\infty$,
one obtains the following well-known Carleman's inequality
\cite{Carlman}, which asserts that for convergent infinite series
$\sum a_n$ with non-negative terms, one has
\begin{equation*}
   \sum^\infty_{n=1}(\prod^n_{k=1}a_k)^{\frac 1{n}}
\leq e\sum^\infty_{n=1}a_n,
\end{equation*}
   with the constant $e$ being best possible.

   It is then natural to study the following weighted version of Carleman's inequality:
\begin{equation}
\label{1}
   \sum^\infty_{n=1}\Big( \prod^n_{k=1}a^{\lambda_k/\Lambda_n}_k \Big )
\leq E\sum^\infty_{n=1}a_n,
\end{equation}
  where the notations are as in \eqref{021}. The task here is to determine the best constant $E$ so that inequality \eqref{1} holds for any convergent infinite series $\sum a_n$ with non-negative terms.  Note that Cartlidge's result (Theorem \ref{thm02}) implies that when \eqref{022} is satisfied, then for any ${\bf a} \in l^p$, one has
\begin{equation}
\label{4}
   \sum^{\infty}_{n=1}\Big{|}\sum^{n}_{k=1}\frac {\lambda_ka_k}{\Lambda_n}
   \Big{|}^p \leq \Big ( \frac {p}{p-L} \Big )^p \sum^{\infty}_{n=1}|a_n|^p.
\end{equation}
  Similar to our discussions above, by a change of variables $a_k \rightarrow a^{1/p}_k$ in \eqref{4} and on letting $p \rightarrow +\infty$, one obtains inequality \eqref{1} with $E=e^{L}$ as long as \eqref{022} is satisfied with $p$ replaced by $+\infty$ there.

  Note that \eqref{022} can be regarded as the case
$p \rightarrow 1^+$ of \eqref{024} while the case $p \rightarrow
+\infty$ of \eqref{024} suggests the following result:
\begin{cor}
\label{thm1}
  Suppose that
\begin{equation*}
  M=\sup_n\frac
    {\Lambda_n}{\lambda_n}\log \Big(\frac {\Lambda_{n+1}/\lambda_{n+1}}{\Lambda_n/\lambda_n} \Big ) < +\infty,
\end{equation*}
  then inequality \eqref{1} holds with $E=e^M$.
\end{cor}
   In fact, the above corollary is a consequence of the following nice result of Bennett (see the proof of \cite[Theorem
   13]{Be1}):
\begin{theorem}
\label{thm6.1}
  Inequality \eqref{1} holds with
\begin{equation*}
  E=\sup_n\frac
    {\Lambda_{n+1}}{\lambda_{n+1}} \prod^n_{k=1} \Big(\frac {\lambda_{k}}{\Lambda_{k}} \Big)^{\lambda_k/\Lambda_n}.
\end{equation*}
\end{theorem}
   It is shown in the proof of Theorem 13 in \cite{Be1} that Corollary \ref{thm1} follows from the
   above theorem. It is also easy to see that $M \leq L$ for $L$ defined by \eqref{022} so that Corollary \ref{thm1} provides a better result
   than what one can infer from Cartlidge's result as discussed above.

   Note that the bound given in Theorem \ref{thm6.1} is global in
   the sense that it involves all the $\lambda_n$'s and it implies
   the local version Corollary \ref{thm1}, in which only the terms
   $\Lambda_n/\lambda_n$ and $\Lambda_{n+1}/\lambda_{n+1}$ are
   involved. It is then natural to ask whether one can obtain a
   similar result for the $l^p$ norms for $p>1$ so that it implies
   the local version Theorem \ref{thm03}. It is our goal in this
   note to present one such result and as our result is motivated
   by the result of Bennett, we will first study the limiting case
   $p \rightarrow +\infty$, namely weighted Carleman's inequalities in the next section before
   we move on to the $l^p$ cases in Section \ref{sec 3}.
\section{A Discussion on Weighted Carleman's Inequalities }
\label{sec 2} \setcounter{equation}{0}
   In this section we study weighted Carleman's inequalities. Our
   goal is to give a different proof of Theorem \ref{thm6.1} than that given in \cite{Be1} and discuss
   some variations of it. It suffices to consider the cases of \eqref{1} with the infinite summations replaced by any finite summations, say from $1$ to $N \geq 1$ here.
   Our starting point is the following result of Pe\v cari\'c and Stolarsky
\cite[(2.4)]{P&S}, which is an outgrowth of Redheffer's approach
in \cite{R1}:
\begin{equation}
\label{5.1}
   \sum^N_{n=1}\Lambda_n(b_n-1)G_n+\Lambda_NG_N
\leq \sum^N_{n=1}\lambda_na_nb^{\Lambda_n/\lambda_n}_n,
\end{equation}
  where ${\bf b}$ is any positive sequence and
\begin{equation*}
   G_n=\prod^n_{k=1}a^{\lambda_k/\Lambda_n}_k.
\end{equation*}
   We now discard the last term on the left-hand side of \eqref{5.1} and make a change of variables $\lambda_na_nb^{\Lambda_n/\lambda_n}_n \mapsto a_n$ to recast
  inequality \eqref{5.1} as
\begin{equation*}
   \sum^N_{n=1}\Lambda_n(b_n-1) \Big( \prod^n_{k=1}\lambda^{-\lambda_k/\Lambda_n}_k \Big )\Big( \prod^n_{k=1}b^{-\Lambda_k/\Lambda_n}_k \Big
   )G_n
\leq \sum^N_{n=1}a_n.
\end{equation*}

   Now, a further change of variables $b_n \mapsto
   \lambda_{n+1}b_n/\lambda_n$ allows us to recast the above
   inequality as
\begin{equation}
\label{5.2}
   \sum^N_{n=1}\Lambda_n\Big(\frac {b_n}{\lambda_n}-\frac {1}{\lambda_{n+1}}\Big) \prod^n_{k=1}b^{-\Lambda_k/\Lambda_n}_k G_n
\leq \sum^N_{n=1}a_n.
\end{equation}

   If one now chooses $b_n=(\Lambda_{n+1}/\lambda_{n+1})/(\Lambda_{n}/\lambda_{n})$ such that
   $\Lambda_n(b_n/\lambda_n-1/\lambda_{n+1})=1$, then one gets
   immediately the following:
\begin{equation*}
   \sum^N_{n=1}\frac {\lambda_{n+1}}{\Lambda_{n+1}} \prod^n_{k=1} \Big( \frac {\Lambda_{k}}{\lambda_k} \Big )^{\lambda_k/\Lambda_n}
   G_n
\leq \sum^N_{n=1}a_n,
\end{equation*}
   from which the assertion of Theorem \ref{thm6.1} can be readily deduced.

   Another natural choice for the values of $b_n$'s is to set
   $\prod^n_{k=1}b^{-\Lambda_k/\Lambda_n}_k=e^{-M}$. From this we
   see that $b_n=e^{M\lambda_n/\Lambda_n}$ and substituting these
   values for $b_n$'s we obtain via \eqref{5.2}:
\begin{equation}
\label{5.2'}
   \sum^N_{n=1}\Lambda_n\Big(\frac {e^{M\lambda_n/\Lambda_n}}{\lambda_n}-\frac {1}{\lambda_{n+1}}\Big) G_n
\leq e^M\sum^N_{n=1}a_n.
\end{equation}
   One checks easily that the above inequality implies Corollary
   \ref{thm1}.

    We now consider a third choice for the $b_n$'s by
    setting
    $b_n=e^{(\Lambda_{n+1}/\lambda_{n+1}-\Lambda_{n}/\lambda_{n})/(\Lambda_{n}/\lambda_{n})}$
    and it follows from this and \eqref{5.2} that
\begin{equation*}
   \sum^N_{n=1}\Lambda_n\Big(\frac {e^{(\Lambda_{n+1}/\lambda_{n+1}-\Lambda_{n}/\lambda_{n})/(\Lambda_{n}/\lambda_{n})}}{\lambda_n}-\frac {1}{\lambda_{n+1}}\Big)
   e^{-\sum^n_{k=1}\frac {\lambda_k}{\Lambda_n}(\frac {\Lambda_{k+1}}{\lambda_{k+1}}-\frac {\Lambda_{k}}{\lambda_{k}})}G_n \leq
\sum^N_{n=1}a_n,
\end{equation*}
   from which we deduce the following
\begin{cor}
\label{cor2}
  Suppose that
\begin{equation*}
  M=\sup_n \sum^n_{k=1}\frac {\lambda_k}{\Lambda_n}(\frac {\Lambda_{k+1}}{\lambda_{k+1}}-\frac {\Lambda_{k}}{\lambda_{k}}) < +\infty,
\end{equation*}
  then inequality \eqref{1} holds with $E=e^M$.
\end{cor}
   We point out here that the above corollary also provides a better result
   than what one can infer from Cartlidge's result and it also follows from
   Theorem \ref{thm6.1}.

   We note here the optimal choice for the $b_n$'s will be to choose them to satisfy
\begin{equation*}
  \Lambda_n\Big(\frac {b_n}{\lambda_n}-\frac {1}{\lambda_{n+1}}\Big)
  \prod^n_{k=1}b^{-\Lambda_k/\Lambda_n}_k=e^{-L}.
\end{equation*}
   In general it is difficult to solve for the $b_n$'s from the
   above equations. But we can solve $b_1$ to get
   $b_1=e^L\lambda_1/((e^L-1)\lambda_2)$ and if we set
   $b_n=(\Lambda_{n+1}/\lambda_{n+1})/(\Lambda_{n}/\lambda_{n})$
   for $n \geq 2$, we can then deduce from \eqref{5.2} that
\begin{equation*}
   e^{-L}G_1+\sum^N_{n=2}\Big ( \frac {\Lambda_2(e^{L}-1)}{\lambda_1e^{L}} \Big )^{\lambda_1/\Lambda_N}\frac {\lambda_{n+1}}{\Lambda_{n+1}} \prod^n_{k=1} \Big( \frac {\Lambda_{k}}{\lambda_k} \Big )^{\lambda_k/\Lambda_n}
   G_n
\leq \sum^N_{n=1}a_n.
\end{equation*}
  Note that this gives an improvement upon Theorem \ref{thm6.1} as
  long as $\lambda_2/\Lambda_2>e^{-L}$. Similarly, one obtains
\begin{equation*}
   G_1+\sum^N_{n=2}\Big ( \frac {\lambda_2(e^{L}-1)}{\lambda_1} \Big )^{\lambda_1/\Lambda_N}\Lambda_n\Big(\frac {e^{M\lambda_n/\Lambda_n}}{\lambda_n}-\frac {1}{\lambda_{n+1}}\Big) G_n
\leq e^M\sum^N_{n=1}a_n.
\end{equation*}

\section{The $l^p$ Cases}
\label{sec 3} \setcounter{equation}{0}
   We now return to the discussions on the general $l^p$ cases. Again it suffices to consider the cases of \eqref{4} with the infinite summations replaced by any finite summations, say from $1$ to $N \geq 1$ here. We
   may also assume that $a_n \geq 0$ for all $n$. With the discussions of the previous section in mind,
   here we seek for an $l^p$ version of \eqref{5.2}. Fortunately this is available by noting that it follows from inequality (4.3) of \cite{G6}
   that
\begin{equation}
\label{5.3}
    \sum_{n=1}^{N}
    \Big(\sum^n_{k=1}w_k \Big )^{-(p-1)}\Big( \frac {w_n^{p-1}}{\lambda^p_n}-\frac {w_{n+1}^{p-1}}{\lambda^p_{n+1}} \Big )
    \Lambda^p_n A_n^{p} \leq \sum_{n=1}^N  a_n^{p},
\end{equation}
  where $w_n$'s are positive parameters and
\begin{equation*}
     A_n=\frac {\sum^n_{k=1}\lambda_ka_k}{\Lambda_n}.
\end{equation*}
   By a change of variables $w_n \mapsto \lambda_nw^{1/(p-1)}_n$,
   we can recast inequality \eqref{5.3} as
\begin{equation*}
    \sum_{n=1}^{N}
    \Big( \frac {\sum^n_{k=1}\lambda_kw^{1/(p-1)}_k}{\Lambda_n} \Big )^{-(p-1)}\Big( \frac {w_n}{\lambda_n}-\frac {w_{n+1}}{\lambda_{n+1}} \Big )
    \Lambda_n A_n^{p} \leq \sum_{n=1}^N  a_n^{p},
\end{equation*}
   With another change of variables, $w_n/w_{n+1} \mapsto b_n$,
   we can further recast the above inequality as
\begin{equation}
\label{5.4}
    \sum_{n=1}^{N}
    \Big( \frac {\sum^n_{k=1}\lambda_k\prod^{n}_{i=k}b^{1/(p-1)}_i}{\Lambda_n} \Big )^{-(p-1)}\Big( \frac {b_n}{\lambda_n}-\frac {1}{\lambda_{n+1}} \Big )
    \Lambda_n A_n^{p} \leq \sum_{n=1}^N  a_n^{p}.
\end{equation}
    Note that if one makes a change of variables $a^p_n \mapsto
   a_n$, then inequality \eqref{5.2} follows from the above
   inequality upon letting $p \rightarrow +\infty$.

   One can then deduce Theorem \ref{thm03} by choosing $b_n=(1-L\lambda_n/(p\Lambda_n))^{-(p-1)}$ in
   \eqref{5.4} (see \cite{G5}). It is easy to see that this gives
   back inequality \eqref{5.2'} upon letting $p \rightarrow +\infty$ by a change of variables $a^p_n
   \mapsto a_n$ and setting $L=M$. We note here that
   the $b_n$'s are so chosen so that the following relations are
   satisfied:
\begin{equation}
\label{3.3}
     \frac
     {\sum^n_{k=1}\lambda_k\prod^{n}_{i=k}b^{1/(p-1)}_i}{\Lambda_n}=
     \frac {p}{p-L}.
\end{equation}

   Now, to get the $l^p$ analogues of Theorem \ref{thm6.1}, we
   just need to note that in the $p \rightarrow +\infty$ case,
   one obtains Corollary \ref{thm1} by setting
   $\prod^n_{k=1}b^{-\Lambda_k/\Lambda_n}_k=e^{-M}$ and the
   conclusion of Corollary \ref{thm1} follows by requiring that $\Lambda_n(b_n/\lambda_n-1/\lambda_{n+1}) \geq
   1$ for the so chosen $b_n$'s. If one instead chooses the
   $b_n$'s so that the conditions $\Lambda_n(b_n/\lambda_n-1/\lambda_{n+1})
   =1$ are satisfied, then Theorem \ref{thm6.1} will follow. Now, in the $l^p$ cases, the choice of
   the $b_n$'s so that the conditions \eqref{3.3} are satisfied implies Theorem \ref{thm03} as
   \eqref{024} implies that for the so chosen $b_n$'s,
\begin{equation*}
     \Big( \frac {b_n}{\lambda_n}-\frac {1}{\lambda_{n+1}} \Big )
    \Lambda_n \geq 1-\frac {L}{p}.
\end{equation*}

   Thus, in order to obtain an result analogue to Theorem \ref{thm6.1} for the $l^p$ cases, we are then motivated to take the $b_n$'s
   so that the following conditions are satisfied:
\begin{equation}
\label{3.4}
     \Big( \frac {b_n}{\lambda_n}-\frac {1}{\lambda_{n+1}} \Big )
    \Lambda_n=1-\frac {L}{p}.
\end{equation}
    We then easily deduce the following $l^p$ analogue of Theorem
    \ref{thm6.1}:
\begin{theorem}
\label{thm3.1}
  Let $1<p<\infty$ be fixed. Let $A$ be a weighted mean matrix given by
    \eqref{021}. If for any integer $n \geq 1$, there exists a positive constant
    $0<L<p$ such that
\begin{equation*}
    \sum^n_{k=1}\frac {\lambda_k}{\Lambda_n}\prod^{n}_{i=k}\Big (\frac {\Lambda_{
    i+1}/\lambda_{i+1}-L/p}{\Lambda_i/\lambda_i} \Big
    )^{1/(p-1)} \leq
     \frac {p}{p-L},
\end{equation*}
    then
    $||A||_{p,p} \leq p/(p-L)$.
\end{theorem}

   It is easy to see by induction that Theorem \ref{thm3.1}
   implies Theorem \ref{thm03}. Of course one should really choose
   the $b_n$'s so that the following relations are satisfied:
\begin{equation*}
    \Big( \frac {\sum^n_{k=1}\lambda_k\prod^{n}_{i=k}b^{1/(p-1)}_i}{\Lambda_n} \Big )^{-(p-1)}\Big( \frac {b_n}{\lambda_n}-\frac {1}{\lambda_{n+1}} \Big )
\Lambda_n=\Big ( \frac {p}{p-L} \Big )^{-p}.
\end{equation*}
   In general it is difficult to determine the $b_n$'s this
   way but one can certainly solve for $b_1$ and by choosing other
   $b_n$'s so that \eqref{3.4} are satisfied, one can obtain a
   slightly better result than Theorem \ref{thm3.1}, we shall
   leave the details to the reader.

   We note here that the choice of the $b_n$'s satisfying \eqref{3.4} corresponds to the following
   choice for the $a_n$'s in Section 4 of \cite{G5} (these $a_n$'s are not to be confused with the $a_n$'s used in the rest of the paper):
\begin{equation*}
  a_n=\Big ( \frac {\Lambda_{n+1}/\lambda_{n+1}-L/p}{\Lambda_n/\lambda_n}  \Big
)^{1/(p-1)}a_{n+1}, \hspace{0.1in} a_1=1.
\end{equation*}

   We also note here that in the case of $\lambda_n=L=1$, on choosing
   $b_n$'s to satisfy \eqref{3.4}, we obtain via \eqref{5.4} that
\begin{equation*}
    \sum_{n=1}^{\infty}
    \Big( \frac {\sum^n_{k=1}\prod^{n}_{i=k}(1+(1-1/p)/i)^{1/(p-1)}}{n} \Big )^{-(p-1)} \Big{(}\frac {1}{n}
\sum^n_{k=1}a_k\Big{)}^p \leq \Big (\frac {p}{p-1} \Big
)\sum_{n=1}^{\infty} a_n^{p}.
\end{equation*}
   This gives an improvement of Hardy's inequality \eqref{eq:1}.

\end{document}